\documentclass[10pt,twoside]{article}
\usepackage{a4,amsmath,amssymb,amsfonts}
\usepackage{Latex-document}

\begin{document}

\bibliographystyle{plain}

 \newcommand{\kkk}[1]{{\Large\bf#1}}
 \newcommand{\nlabel}[1]{\label{#1}\kkk{#1}}

 \def\beginProof{\par \noindent {\bf Proof }}
 \def\endProof{${\bf Q.E.D.}$\par}
 \def\ar#1{\widehat{#1}}
 \def\mn{{\mu_{n}}}
 \def\pr{^{\prime}}
 \def\prpr{^{\prime\prime}}
 \def\mtr#1{\overline{#1}}
 \def\ra{\rightarrow}
 \def\Bbb{\bf}
 \def\mC{{\mathbb C}}
 \def\mF{{\mathbb F}}
 \def\mG{{\mathbb G}}
 \def\mP{{\mathbb P}}
 \def\mQ{{\mathbb Q}}
 \def\mR{{\mathbb R}}
 \def\mZ{{\mathbb Z}}
 \def\mC{{\mathbb C}}
 \def\mN{{\mathbb N}}
 \def\Qb{\mtr{\mathbb Q}}
 \def\Zn{{{\mathbb Z}/n}}
 \def\Qmn{{\mQ(\mn)}}
 \def\refeq#1{(\ref{#1})}
 \def\umn{^{\mn}}
 \def\lmn{_{\mn}}
 \def\blb{{\big(}}
 \def\brb{{\big)}}
 \def\Aut{\mathop{\rm Aut}\nolimits}
 \def\Hom{\mathop{\rm Hom}\nolimits}
 \def\Tr{\mathop{\rm Tr}\nolimits}
 \def\Gal{\mathop{\rm Gal}\nolimits}
 \def\End{\mathop{\rm End}\nolimits}
 \def\Det{\mathop{\rm Det}\nolimits}
 \def\Proj{\mathop{\rm Proj}\nolimits}
 \def\Spec{\mathop{\rm Spec}\nolimits}
 \def\deg{\mathop{\rm deg}\nolimits}
 \def\mod{\mathop{\rm mod}\nolimits}
\def\Ker{\mathop{\rm Ker}\nolimits}
\def\Tor{\mathop{\rm Tor}\nolimits}
\def\Stab{\mathop{\rm Stab}\nolimits}
\def\Lie{\mathop{\rm Lie}\nolimits}
\def\codim{\mathop{\rm codim}\nolimits}
\def\proj{\mathop{\rm pr}\nolimits}
\def\PGL{\mathop{\rm PGL}\nolimits}
\def\et{{\rm et}}
\def\reg{{\rm reg}}
\def\cl{{\rm cl}}
\def\id{{\rm id}}
\let\phi\varphi
\let\setminus\smallsetminus
\let\into\hookrightarrow
\def\Fp{{\mathfrak p}}

 \newtheorem{theor}{Theorem}[section]
 \newtheorem{prop}[theor]{Proposition}
 \newtheorem{cor}[theor]{Corollary}
 \newtheorem{lemma}[theor]{Lemma}
 \newtheorem{sublem}[theor]{sublemma}
 \newtheorem{defin}[theor]{Definition}
 \newtheorem{conj}[theor]{Conjecture}
 \newtheorem{rem}[theor]{Remark}
\renewcommand{\thesection}{\arabic{section}}
\def\Section#1{\section{\hskip -1em . \hskip 0.6em #1}}
\def\volumeno{I}
\markboth{On Hrushovski's Proof of the Manin-Mumford Conjecture}{Richard
Pink \quad Damian Roessler}
\title{\bf On Hrushovski's Proof of the \vskip -2mm Manin-Mumford Conjecture\vskip 6mm}
\author{{\bf Richard Pink}\thanks{Department of Mathematics,
ETH-Zentrum, CH-8092 Z\"urich, Switzerland. E-mail: pink@math.ethz.ch}
\quad Damian Roessler\thanks{Department of Mathematics, ETH-Zentrum,
CH-8092 Z\"urich, Switzerland. E-mail:
roessler@math.ethz.ch}\vspace*{-0.5cm}}
\date{\vspace{-8mm}}

\maketitle

\thispagestyle{first} \setcounter{page}{539}

\begin{abstract}

\vskip 3mm

The Manin-Mumford conjecture in characteristic zero was first proved
by Raynaud.  Later, Hrushovski gave a different proof using model
theory.  His main result from model theory, when applied to abelian
varieties, can be rephrased in terms of algebraic geometry.  In this
paper we prove that intervening result using classical algebraic
geometry alone.  Altogether, this yields a new proof of the
Manin-Mumford conjecture using only classical algebraic geometry.

\vskip 4.5mm

\noindent {\bf 2000 Mathematics Subject Classification:} 14K12.

\noindent {\bf Keywords and Phrases:} Torsion points, Abelian varieties,
Manin-Mumford conjecture.
\end{abstract}

\vskip 12mm

\Section{Introduction}

\vskip-5mm \hspace{5mm}

The Manin-Mumford conjecture is the following statement.

\begin{theor}
\label{MM}
Let $A$ be an abelian variety defined over~$\Qb$ and $X$ a closed
subvariety of~$A$.  Denote by $\Tor(A)$ the set of torsion points
of~$A$. Then
$$X\cap \Tor(A) = \bigcup_{i\in I}X_i\cap \Tor(A),$$
where $I$ is a finite set and each $X_i$ is the translate by an
element of $A$ of an abelian subvariety of~$A$, immersed in~$X$.
\end{theor}

This conjecture has a long history of proofs.  A first partial proof
was given by Bogomolov in \cite{Bo}, who proved the statement when
$\Tor(A)$ is replaced by its $\ell$-primary part for a prime
number~$\ell$; he applies results of Serre, Tate and Raynaud (see
\cite{Tate}) on the existence of Hodge-Tate module structures on the
Tate module of abelian varieties over discrete valuation rings.
A full proof of the conjecture was then given by Raynaud in \cite{Ra2}
(see \cite{Ra1} for the case $\dim(X)=1$);
his proof follows from a study of the
reduction of $A$ modulo $p^2$.
A third and full proof of the conjecture was given by Hrushovski in
\cite{HR}; he uses Weil's theorem on the characteristic polynomial of
the Frobenius morphism on abelian varieties over finite fields in
conjunction with the model theory (of mathematical logic) of the
theory of fields with a distinguished automorphism.
A fourth proof of the conjecture was given by Ullmo and Zhang (based
on ideas of Szpiro) in \cite{Ullmo} and \cite{Zhang} and goes via
diophantine approximation and Arakelov theory.  They actually prove a
more general conjecture of Bogomolov which generalizes the
Manin-Mumford conjecture.
\medskip

This article was inspired by Hrushovski's proof.  The bulk of
Hrushovski's proof lies in the model theory part.  It culminates in a
result which, when applied to the special case of abelian varieties
and stripped of model theoretic terminology, is essentially
Theorem~\ref{HMain} below.  In Section~2 we prove Theorem~\ref{HMain}
with classical algebraic geometry alone.  Neither scheme theory, nor
Arakelov theory, nor mathematical logic are used.  In section~3, for
the sake of completeness, we show how to apply Theorem~\ref{HMain} to
prove the Manin-Mumford conjecture.  \medskip

In a subsequent article (cf. \cite{PR2}), we shall consider
the analogue of \ref{HMain} for varieties over function fields
of characteristic $p>0$.


\Section{Hrushovski's theorem for abelian varieties}

\vskip-5mm \hspace{5mm}

Let $K$ be an algebraically closed field of characteristic zero, endowed
with an automorphism~$\sigma$.  Let $A$ be an abelian variety over~$K$
and $X$ a closed subvariety of~$A$.  For ease of notation, we use the
language of classical algebraic geometry; thus $A$ and $X$ denote the
respective sets of $K$-valued points.  We assume that $X\subset A$ are
defined already over the fixed field~$K^\sigma$.  The automorphism of
$A$ induced by $\sigma$ is again denoted by~$\sigma$. Let
$P(T)\in\mZ[T]$ be a monic polynomial with integral coefficients. In
\cite[Cor.~4.1.13, p.90]{HR}, Hrushovski proves the generalisation of
the following theorem to semi-abelian varieties.

\begin{theor}
Let $\Gamma$ denote the kernel of the homomorphism $P(\sigma): A\ra
A$.  Assume that no complex root of $P$ is root of unity.  Then
$$X\cap\Gamma = \bigcup_{i\in I}X_i\cap\Gamma,$$
where $I$ is a finite set and each $X_i$ is the translate by an
element of $A$ of an abelian subvariety of~$A$, immersed in~$X$.
\label{HMain}
\end{theor}

\noindent {\bf Remark} If roots of unity are not excluded, the group
$\Gamma$ becomes too large for such a result.  For example, if $T^m-1$
divides~$P$, all points of $A$ over the fixed field $K^{\sigma^m}$ of
$\sigma^m$ are contained in~$\Gamma$.
\medskip

\beginProof \,\,Write $P(T) = \sum_{i=0}^n a_i T^i$ with
$a_i\in\mZ$ and $a_n=1$.  Let $F$ be the endomorphism of $A^n$ defined
by the matrix
\begin{displaymath}
\left(
\begin{array}{ccccc}
0 & 1 & 0 & \ldots & 0 \cr
\vdots & \ddots & \ddots & \ddots & \vdots \cr
\vdots &  & \ddots & \ddots & 0 \cr
0 & \ldots & \ldots & 0 & 1 \cr
-a_0 & -a_1 & \ldots & \ldots & -a_{n-1}
\end{array}
\right)
\end{displaymath}
which is the companion matrix of the polynomial~$P$, and note that
$P(F)=0$.  In the obvious way $\sigma$ induces an automorphism
of~$A^n$, denoted again by~$\sigma$, that sends $X^n$ to itself.  Let
$\Delta$ denote the kernel of the homomorphism $F-\sigma: A^n\ra A^n$.
By construction, there is a canonical bijection
$$\Gamma\to\Delta,\ x\mapsto (x,\sigma(x),\dots,\sigma^{n-1}(x)).$$
Since $\sigma(X)=X$, this induces a bijection
$$X\cap\Gamma \to X^n\cap\Delta.$$
Its inverse is given by the projection to the first factor $A^n\to A$.
Clearly, we are now reduced to the following theorem, applied to
$X^n\subset A^n$ in place of $X\subset A$.

\begin{theor}
Let $F:A\to A$ be an algebraic endomorphism that commutes with
$\sigma$ and that satisfies $P(F)=0$.  Let $\Delta$ denote the kernel
of the homomorphism $F-\sigma: A\ra A$.  Assume that no complex root
of $P$ is root of unity.  Then
$$X\cap\Delta = \bigcup_{i\in I}X_i\cap\Delta,$$
where $I$ is a finite set and each $X_i$ is the translate by an
element of $A$ of an abelian subvariety of~$A$, immersed in~$X$.
\label{HMain2}
\end{theor}

\begin{rem}
\rm If $K$ is algebraic over the fixed field~$K^\sigma$, every element
$a\in\Delta$ satisfies $F^m(a) = \sigma^m(a) = a$ for some $m\ge1$.
In other words, we have $a \in \Ker(F^m-\id)$.  The assumptions on $F$
and $P$ imply that $F^m-\id$ is an isogeny; hence $a$ is a torsion
element.  It follows that $\Delta$ is a torsion subgroup, and the
theorem follows from the Manin-Mumford conjecture in this case.
However, the scope of the above theorem is somewhat wider, and the
Manin-Mumford conjecture will be deduced from it.
\end{rem}

\beginProof
\,\,Let $Y$ be the Zariski closure of $X\cap\Delta$.  We claim
that $\sigma(Y) = F(Y) = Y$.  To see this, note first that
$\sigma$ commutes with~$F$, and so $\sigma(\Delta) = \Delta$.  By
assumption we have $\sigma(X)=X$; hence $\sigma(X\cap\Delta) =
X\cap\Delta$.  Since $\sigma:A\to A$ comes from an automorphism of
the underlying field, it is a homeomorphism for the Zariski
topology, so we have $\sigma(Y)=Y$. On the other hand, the maps
$\sigma$ and $F$ coincide on~$\Delta$, which implies
$F(X\cap\Delta) = X\cap\Delta$.  As $F$ is a proper algebraic
morphism, we deduce that $F(Y)=Y$. Clearly, Theorem~\ref{HMain2}
is now reduced to the following theorem (see \cite[Th. 3]{Bogo}
for the case where $F$ is the multiplication by an integer $n>1$).

\begin{theor}
Let $A$ be an abelian variety over an algebraically closed field of
characteristic zero.  Let $F:A\to A$ be an algebraic endomorphism none
of whose eigenvalues on $\Lie A$ is a root of unity.  Let $Y$ be a
closed subvariety of $A$ satisfying $F(Y)=Y$.  Then $Y$ is a finite
union of translates of abelian subvarieties of~$A$.
\label{HMain3}
\end{theor}

\beginProof
\,\,We proceed by induction on the dimension $d$ of~$A$.  For
$d=0$ the statement is obvious; hence we assume $d>0$.  Next
observe that $Y\subset F(A)$.  Thus if $F$ is not surjective, we
can replace $A$ by $F(A)$ and $Y$ by $F(A)\cap Y$, and are
finished by induction.  Thus we may assume that $F$ is an isogeny.
\medskip

Since $F$ is proper and $F(Y)=Y$, every irreducible component of $Y$
is the image under $F$ of some irreducible component of~$Y$.  Since
the set of these irreducible components is finite, it is therefore
permuted by~$F$.  We fix such an irreducible component $Z$ and an
integer $m\ge1$ such that $F^m(Z)=Z$.
\medskip

Any power $F^{rm}:A\to A$ is an isogeny, and since $\mathop{\rm
char}(K)=0$, it is a separable isogeny.  As a morphism of schemes it
is therefore a finite \'etale Galois covering with Galois group
$\Ker(F^{rm})$, acting by translations on~$A$.  The same follows for
the induced covering $(F^{rm})^{-1}(Z) \to Z$.  As $Z$ is irreducible,
the irreducible components of $(F^{rm})^{-1}(Z)$ are transitively
permuted by $\Ker(F^{rm})$ and each of them has dimension $\dim(Z)$.
Since $F^{rm}(Z)=Z$, we have $Z\subset (F^{rm})^{-1}(Z)$, and so $Z$
itself is one of these irreducible components.  Let $G_r$ denote the
stabilizer of $Z$ in $\Ker(F^{rm})$.  Then $F^{rm}:Z\to Z$ is
generically a finite \'etale Galois covering with Galois group~$G_r$.
\medskip

We now distinguish two cases.  Let $\Stab_A(Z)$ denote the translation
stabilizer of $Z$ in~$A$, which is a closed algebraic subgroup of~$A$.

\begin{lemma}
If $|G_1|>1$, then $\dim(\Stab_A(Z))>0$.
\end{lemma}

\beginProof
\,\,For any $r\ge1$, the morphism $F^{rm}:Z\to Z$ is finite
separable of degree $|G_r|$.  Since degrees are multiplicative in
composites, this degree is also equal to $|G_1|^r$.  Thus if
$|G_1|>1$, we find that $|G_r|$ becomes arbitrarily large
with~$r$.  Therefore $\Stab_A(Z)$ contains arbitrarily large
finite subgroups, so it cannot be finite.
\endProof

\begin{lemma}
If $|G_1|=1$, then $\dim(\Stab_A(Z))>0$ or $Z$ is a single point.
\label{lemma_aut}
\end{lemma}

\beginProof
\,\,Since $Z$ is irreducible, the assertion is obvious when
$\dim(Z)=0$. So we assume that $\dim(Z)>0$.  The assumption
$|G_1|=1$ implies that $F^m$ induces a finite (separable) morphism
$\phi: Z\to Z$ of degree $1$.  The set of fixed points of any
positive power $\phi^r$ is then $Z\cap\Ker(F^{rm}-\id)$.  On the
other hand the assumptions on $F$ and $P$ imply that $F^{rm}-\id$
is an isogeny on~$A$.  Thus this fixed point set is finite for
every $r\ge1$; hence $\phi$ has infinite order.
\smallskip

Assume now that $\Stab_A(Z)$ is finite.  By Ueno's theorem
\cite[Thm.~3.10]{Ueno} $Z$ is then of general type, in the sense that
any smooth projective variety birationally equivalent to $Z$ is of
general type (see \cite[Def.~1.7]{Ueno}).  But the group of birational
automorphisms of any irreducible projective variety $Z$ of general
type is finite (see \cite{Matsu} and also \cite[Th.  10.11]{II}) which
yields a contradiction.
\endProof
\medskip

If $Z$ is a single point, it is a translate of the trivial abelian
subvariety of~$A$.  Otherwise we know from the lemmas above that
$\dim(\Stab_A(Z))>0$, and we can use the induction hypothesis.  Let
$B$ denote the identity component of $\Stab_A(Z)$.  Since $F^m(Z) =
Z$, we also have $F^m(B) = B$.  Set $\bar A:=A/B$ and $\bar Z:=Z/B$,
and let $\bar F$ denote the endomorphism of $\bar A$ induced by~$F^m$.
Then we have $\bar F(\bar Z)=\bar Z$, and every eigenvalue of $\bar F$
on $\Lie\bar A$ is an eigenvalue of $F^m$ on $\Lie A$ and therefore
not a root of unity.  By Theorem~\ref{HMain3}, applied to $(\bar
A,\bar Z,\bar F)$ in place of $(A,Y,F)$, we now deduce that $\bar Z$
is a finite union of translates of abelian subvarieties of~$\bar A$.
But $Z$ is irreducible; hence so is~$\bar Z$.  Thus $\bar Z$ itself is
a translate of an abelian subvariety.  That abelian subvariety is
equal to $A'/B$ for some abelian subvariety $A'\subset A$
containing~$B$, and so $Z$ is a translate of~$A'$, as desired.  This
finishes the proof of Theorem~\ref{HMain3}, and thus also of Theorems
\ref{HMain2} and~\ref{HMain}.
\endProof
\medskip

\noindent {\bf Remark} Suppose that in the statement of the preceding
theorem, we replace the assumption that none of the eigenvalues of $F$
on $\Lie A$ is a root of unity by the weaker assumption that none of the
eigenvalues of $F$ on $\Lie A$ is an algebraic unit.  This case is
sufficient for the application to the Manin-Mumford conjecture.  Under
this weaker assumption, the following alternative proof of Lemma
\ref{lemma_aut} can be given; it does not use the theorems of Ueno and
Matsumura but only elementary properties of cycle classes.
\smallskip

We work under the hypotheses of Lemma \ref{lemma_aut} and with the
notations used above.  First, everything can be defined over a
countable algebraically closed subfield of~$K$, and this subfield can
be embedded into~$\mC$; thus without loss of generality we may assume
that $K=\mC$.  Then for every integer $i\ge0$ we abbreviate $H^i :=
H^i(A,\mZ)$.  Let $c:=\codim_A(Z)$ and let $\cl(Z) \in H^{2c}$ be the
cycle class of $Z \subset A$.  We calculate the cycle class of
$(F^{m})^{-1}(Z)$ in two ways.  On the one hand, we have seen that the
group $\Ker(F^m)$ acts transitively on the set of irreducible
components of $(F^{m})^{-1}(Z)$; the assumption $|G_1|=1$ implies that
it also acts faithfully.  Thus the number of irreducible components is
$|\Ker(F^m)|$.  Recall that $Z$ is one of them.  Since translation on
$A$ does not change cycle classes, we find that all irreducible
components of $(F^{m})^{-1}(Z)$ have cycle class $\cl(Z)$; hence
$\cl\bigl((F^{m})^{-1}(Z)\bigr) = |\Ker(F^m)| \cdot \cl(Z)$.  On the
other hand $F$ induces a pullback homomorphism $F^*: H^i\to H^i$ for
every $i\ge0$, and by functoriality of cycle classes we have
$\cl\bigl((F^{m})^{-1}(Z)\bigr) = F^*\bigl(\cl(Z)\bigr)$.  As $\cl(Z)$
is non-zero, we deduce that $|\Ker(F^m)|$ is an eigenvalue of $F^*$ on
$H^{2c}$.
\smallskip

Let $d:=\dim(A)$, which is also the codimension of any point in~$A$.
Repeating the above calculation with the cycle class of a point, we
deduce that $|\Ker(F^m)|$ is also an eigenvalue of $F^*$ on the
highest non-vanishing cohomology group $H^{2d}$.  Since cup product
yields isomorphisms $H^i \cong \Lambda^i H^1$ for all $i\ge0$, which
are compatible with~$F^*$, it yields an $F^*$-equivariant perfect
pairing $H^{2c}\times H^{2(d-c)} \to H^{2d}$.  From this we deduce
that $1$ is an eigenvalue of $F^*$ on $H^{2(d-c)}$.  Moreover, this
eigenvalue must be a product of $2(d-c)$ eigenvalues of $F^*$
on~$H^1$.  Now the eigenvalues of $F^*$ on $H^1$ are precisely those
of $F$ on $\Lie A$ and their complex conjugates.  By assumption, they
are algebraic integers but no algebraic units.  Thus a product of such
numbers can be $1$ only if it is the empty product.  This shows that
$2(d-c)=0$, that is, $\codim_A(Z) = \dim(A)$; hence $Z$ is a point, as
desired.
\endProof

\Section{Proof of the Manin-Mumford conjecture}

\vskip-5mm \hspace{5mm}

As a preparation let $q=p^r$, where $p$ is a prime number and $r\geq
1$.  Let $\mF_q$ be the unique field with $q$ elements.  Let $A$ be an
abelian variety defined over~$\mF_q$.  As in classical algebraic
geometry we identify $A$ with the set of its $\mtr{\mF_q}$-valued
points.  Let $\ell$ be a prime number different from~$p$, and let
$T_\ell(A)$ be the $\ell$-adic Tate module of~$A$, which is a free
$\mZ_\ell$-module of rank $2\dim(A)$.  We denote by $\phi$ the
Frobenius morphism $A\ra A$, which acts on the coordinates of points
by taking $q$-th powers.  It also acts on the Tate module via its
action on the torsion points.  The following result, an analogue of
the Riemann hypothesis, is due to Weil (see \cite{Weil}):

\begin{theor}
Let $T$ be an indeterminate.  The characteristic polynomial of $\phi$
on $T_\ell(A)\otimes_{\mZ_\ell}\mQ_\ell$ is a monic polynomial in
$\mZ[T]$.  It is independent of~$\ell$, and all its complex roots have
absolute value~$\sqrt{q}$.
\label{WeilRH}
\end{theor}

Consider now an abelian variety $A$
defined over~$\Qb$ and a closed subvariety $X$ of~$A$.  We choose a
number field $L\subset\Qb$ over which both $X\subset A$ can be defined
and fix their models over~$L$.  For any abelian group $G$ we write
$\Tor(G)$ for the group of torsion points of~$G$.  Moreover, for any
prime $p$ we write $\Tor_p(G)$ for the subgroup of torsion points of
order a power of~$p$, and $\Tor^p(G)$ for the subgroup of torsion
points of order prime to~$p$.  Note that $\Tor(G) = \Tor^p(G) \oplus
\Tor_p(G)$.
\medskip

Choose a prime ideal $\mathfrak p$ of ${\cal O}_L$ where $A$ has good
reduction.  Let $\mF_q$ be the finite field ${\cal O}_L/{\mathfrak
p}$, where $q$ is a power of a prime number~$p$.
The following lemma is lemma 5.0.10 in \cite[p. 105]{HR}; we reproduce
the proof for the convenience of the reader. We use Weil's
theorem \ref{WeilRH} and reduction modulo $\mathfrak p$ to
obtain an automorphism of $\mtr{\mQ}$ and a polynomial that we
can feed in Theorem~\ref{HMain} to obtain Theorem~\ref{MM}.

\begin{lemma}
\label{Frob}
There is an element $\sigma_{\mathfrak p}\in\Gal(\mtr{\mQ}|L)$ and a
monic polynomial $P_\Fp(T)\in\mZ[T]$ all of whose complex roots have
absolute value~$\sqrt{q}$, such that $P_\Fp(\sigma_{\mathfrak
p})(x)=0$ for every $x\in\Tor^p(A)$.
\end{lemma}

\beginProof
\,\,Let $L_\Fp$ be the completion of $L$ at~$\Fp$.  Extend the
embedding $L\into L_\Fp$ to an embedding of the algebraic closures
$\mtr{\mQ} = \mtr{L} \into \mtr{L_{\mathfrak p}}$, and the
surjection ${\cal O}_{L_\Fp}\ra \mF_q$ to a surjection ${\cal
O}_{\mtr{L_\Fp}} \to \mtr{\mF_q}$.  On the prime-to-$p$ torsion
groups we then obtain natural isomorphisms
$$
\Tor^p(A) \cong
\Tor^p(A_{\mtr{L_{\mathfrak p}}}) \cong
\Tor^p(A_{\mtr{\mF_q}}).
$$
The second isomorphism expresses the fact that the field of definition
of every prime-to-$p$ torsion point is unramified at~$\Fp$.
\medskip

Now, as before, let $\phi$ denote the automorphism of $\mtr{\mF_q}$
and of $\Tor^p(A_{\mtr{\mF_q}})$ induced by the Frobenius morphism
over~$\mF_q$.  Then $\phi$ can be lifted to an element $\sigma_\Fp\in
\Gal(\mtr{\mQ}|L)$ making the above isomorphisms equivariant.  To see
this, one first lifts $\phi$ to an element $\tau_\Fp^{\rm nr}$ of
$\Gal(L^{\rm nr}_{\mathfrak p}| L_{\mathfrak p})$, where $L_{\mathfrak
p}^{\rm nr}$ is the maximal unramified extension of~$L_{\mathfrak p}$.
This lifting exists and is unique by \cite[Th.~1, p.~26]{CaFo}.  This
element can then be lifted further to a (non-unique) element
$\tau_{\mathfrak p}$ of $\Gal(\mtr{L_{\mathfrak p}}|{L_{\mathfrak
p}})$, since $L^{\rm nr}_{\mathfrak p}$ is a subfield of
$\mtr{L_{\mathfrak p}}$.  By construction the action of $\tau_\Fp$ on
$\Tor^p(A_{\mtr{L_{\mathfrak p}}})$ corresponds to the action of
$\phi$ on $\Tor^p(A_{\mtr{\mF_q}})$.  The restriction of $\tau_\Fp$ to
$\mtr{\mQ}$ gives the desired element~$\sigma_\Fp$.
\medskip

Let $P_\Fp(T)$ be the characteristic polynomial of $\phi$ on
$T_\ell(A_{\mtr{\mF_q}})\otimes_{\mZ_\ell}\mQ_\ell$ for any prime
number $\ell\not=p$.  By construction, we then have $P_\Fp(\phi)=0$ on
$\Tor_\ell(A_{\mtr{\mF_q}})$.  By Weil's result quoted above, the same
equation holds for every prime $\ell\not=p$; so it holds on
$\Tor^p(A_{\mtr{\mF_q}})$.  From the construction of $\sigma_\Fp$ we
deduce that $P_\Fp(\sigma_\Fp)=0$ on $\Tor^p(A)$.  Finally, Weil's
result also describes the complex roots of~$P_\Fp$.
\endProof
\medskip

Let now $L^p$, $L_p\subset\Qb$ be the fields generated over $L$ by the
coordinates of all points in $\Tor^p(A)$, resp.\ in $\Tor_p(A)$.  Both
are infinite Galois extensions of~$L$.  Their intersection is known to
be finite over $L$ by Serre \cite[pp.~33--34, 56--59]{Serre}.  Thus
after replacing $L$ by $L^p\cap L_p$, we may assume that $L^p$ and
$L_p$ are linearly disjoint over~$L$.  The subfield of $\Qb$ generated
by the coordinates of all points in $\Tor(A) = \Tor^p(A) \oplus
\Tor_p(A)$ is then canonically isomorphic to $L^p\otimes_LL_p$.
\medskip

Let ${\mathfrak p}\pr$ be a second place of good reduction of~$A$, of
residue characteristic different from~$p$.  Let $\sigma_\Fp$,
$\sigma_{\Fp'}$ and $P_\Fp$, $P_{\Fp'}$ be the automorphisms and
polynomials provided by Lemma~\ref{Frob}, applied to $\Fp$, resp.\
to~$\Fp'$.  The automorphism of $L^p\otimes_LL_p$ induced by
$\sigma_\Fp\otimes \sigma_{\Fp'}$ then extends to some automorphism
$\sigma$ of $\Qb$ over~$L$.  Since $P_\Fp(\sigma_\Fp)$ vanishes on
$\Tor^p(A)$, so does $P_\Fp(\sigma)$.  Similarly,
$P_{\Fp'}(\sigma_{\Fp'})$ vanishes on $\Tor_p(A) \subset
\Tor^{p\pr}(A)$; hence so does $P_{\Fp'}(\sigma)$.  Thus with
$P(T):=P_\Fp(T) P_{\Fp'}(T)$ we deduce that $P(\sigma)$ vanishes on
$\Tor(A)$.  In other words, we have $\Tor(A) \subset \Gamma := \Ker
P(\sigma)$.  With $K=\Qb$ Theorem~\ref{MM} now follows directly from
Theorem~\ref{HMain}.
\endProof


\label{lastpage}

\end{document}